\documentclass[12pt]{amsart}
\usepackage{amsmath,amssymb,latexsym}
\def\dcl{\mathrm{dcl}}

\def\acl{\mathrm{acl}}
\def\bdd{\mathrm{bdd}}
\def\cb{\mathrm{Cb}}

\def\M{\mathfrak M}

\def\tp{\mathrm{tp}}
\def\lstp{\mathrm{Lstp}}
\def\stp{\mathrm{stp}}

\def\stab{\mathrm{stab}}

\def\Ind#1#2{#1\setbox0=\hbox{$#1x$}\kern\wd0\hbox to 0pt{\hss$#1\mid$\hss}
\lower.9\ht0\hbox to 0pt{\hss$#1\smile$\hss}\kern\wd0}
\def\ind{\mathop{\mathpalette\Ind{}}}
\def\Notind#1#2{#1\setbox0=\hbox{$#1x$}\kern\wd0\hbox to 0pt{\mathchardef
\nn="3236\hss$#1\nn$\kern1.4\wd0\hss}\hbox to 0pt{\hss$#1\mid$\hss}\lower.9\ht0
\hbox to 0pt{\hss$#1\smile$\hss}\kern\wd0}
\def\nind{\mathop{\mathpalette\Notind{}}}

\theoremstyle{plain}
\newtheorem{theorem}{Theorem}[section]
\newtheorem{proposition}[theorem]{Proposition}
\newtheorem{fact}[theorem]{Fact}
\newtheorem{lemma}[theorem]{Lemma}
\newtheorem{corollary}[theorem]{Corollary}
\newtheorem*{claim}{Claim}
\newtheorem*{SL}{Socle Lemma}
\newtheorem*{SC}{Corollary}

\theoremstyle{definition}
\newtheorem{definition}[theorem]{Definition}
\newtheorem{remark}[theorem]{Remark}
\newtheorem{expl}[theorem]{Example}
\newtheorem*{question}{Question}

\def\bsp{\begin{expl}}
\def\ebsp{\end{expl}}
\def\beh{\begin{claim}}
\def\ebeh{\end{claim}}
\def\defn{\begin{definition}}
\def\edefn{\end{definition}}
\def\satz{\begin{theorem}}
\def\esatz{\end{theorem}}
\def\tats{\begin{fact}}
\def\etats{\end{fact}}
\def\kor{\begin{corollary}}
\def\ekor{\end{corollary}}
\def\lmm{\begin{lemma}}
\def\elmm{\end{lemma}}
\def\bem{\begin{remark}}
\def\ebem{\end{remark}}
\def\bew{\par\noindent{\em Proof: }}
\def\bewbeh{\par\noindent{\em Proof of Claim: }}
\def\satzli{\begin{proposition}}
\def\esatzli{\end{proposition}}
\def\frag{\begin{question}}
\def\efrag{\end{question}}

\begin{document}
\title{Rigidity, internality and analysability}
\author{Daniel Palac\'\i n}
\author{Frank O. Wagner}

\address{Universit\'e de Lyon; CNRS; Universit\'e Lyon 1; Institut Camille Jordan UMR5208, 43 bd du 11 novembre 1918, 69622 Villeurbanne Cedex, France}
\email{palacin@math.univ-lyon1.fr}
\email{wagner@math.univ-lyon1.fr}
\date{}
\thanks{Partially supported by ANR-09-BLAN-0047 Modig. The first author was also partially supported by research project MTM 2011-26840 of the Spanish government and
research project 2009SGR 00187 of the Catalan government. \newline Both authors thank Anand Pillay for a discussion around the definition of the Canonical Base Property in the general setting.}
\begin{abstract} We prove a version of Hrushovski's {\em Socle Lemma} for rigid groups in an arbitrary simple theory.\end{abstract}
\maketitle

\section{Introduction}

One of the main consequences of the {\em canonical base property} (in short, CBP) for a simple theory of finite $SU$-rank is a certain generalization of the the so-called {\em Socle Lemma}, due originally to Hrushovski for suitably rigid groups (\cite[Proposition 4.3]{udi} and \cite[Proposition 3.6.2]{udi2}). Namely, in a group of finite $SU$-rank with the CBP, every type with finite stabilizer is almost internal to the family of types of rank $1$. This was noted in \cite{pz} and elaborated in \cite{hpp}.

The formulation of the canonical base property in model-theoretic terms was influenced by the results of Campana \cite{ca} and Fujiki \cite{fu} in compact complex spaces and an analogous results due to Pillay and Ziegler \cite{pz} on differential (and difference) algebraic varieties in characteristic $0$. In particular, the group-like version in the case of differentially closed fields of characteristic $0$ yields an account of the Mordell-Lang Conjecture for function fields in characteristic $0$ without using Zariski geometries, see \cite{pz}.
The CBP, named by Moosa and Pillay \cite{mp}, states that for any tuple $a$ of finite $SU$-rank and any $b$, the type of the canonical base of $\tp(a/b)$ over $a$ is almost internal to the family of types of rank $1$. Clearly, this is a property of the finite-rank context; nevertheless, not all theories of finite rank satisfy the CBP \cite{hpp}.

A better understanding of the CBP was provided in \cite{zoe}, where Chatzidakis studied simple theories with the CBP in general. Extending the work of Pillay and Ziegler, she moreover showed that the theory of existentially closed difference fields of any characteristic has the CBP. On the other hand, Chatzidakis proved that replacing internality by analysability in the statement of the CBP, and considering the family of non one-based types of rank $1$, one obtains a weak version of the CBP which is satisfied by all types of finite rank in any simple theory.

Recall that a partial type $\pi$ over a set $A$ in a simple theory is {\em one-based} if for any tuple $\bar a$ of realizations of $\pi$ and any $B\supseteq A$ the canonical base of $\tp(\bar a/B)$ over $A\bar a$ is bounded. One-basedness implies that the forking geometry is
particularly well-behaved; for instance one-based groups are {\em rigid}. That is, every subgroup is commensurable with one hyperdefinable over $\bdd(\emptyset)$ --- see Definition \ref{def:rigid}. It turns out that the CBP and the weak version of the CBP correspond to a relative version of one-basedness with respect to the family of types of rank $1$. This connection was noticed by Kowalski and Pillay in \cite{kp}, and used to describe the structure of type-definable groups in stable theories satisfying the CBP.
The relation between one-basedness and the weak CBP was then used by the authors \cite{pw} to generalize and study the weak CBP outside the finite $SU$-rank setting by replacing the family of types of rank $1$ by an arbitrary family of partial types. These ideas also appear in \cite{bm-pw} where Blossier, Mart\'\i n-Pizarro and the second author study a generalization of the CBP in a different direction.

In \cite{hpp} Hrushovski, Pillay and the first author proved the CBP for non-multidimensional theories of finite Morley rank, assuming that all Galois groups are rigid. The aim of this paper is to obtain directly the Socle Lemma for rigid groups in arbitrary simple theories, without passing through the CBP (and hence without any assumption of finiteness of $SU$-rank). In addition, we remark that the non-multidimensionality assumption of \cite{hpp} is not required.

Our notation is standard and follows \cite{wa00}. Throughout the paper, the ambient theory will be simple, and we shall
be working in $\M^{heq}$, where $\M$ is a sufficiently saturated model of the ambient theory. Thus $A$, $B$, $C$,\ldots will denote hyperimaginary parameter sets, and $a$, $b$, $c$,\ldots tuples of hyperimaginaries. Moreover, $\dcl=\dcl^{heq}$.

\section{The canonical base property}

For the rest of the paper $\Sigma$ will be an $\emptyset$-invariant family of partial types. We shall think of $\Sigma$ both as a set of partial types (as in Definition \ref{defn} below), or as the $\emptyset$-invariant set of its realisations (as in $\acl(\Sigma)$, $\bdd(\Sigma)$ or $a\ind\Sigma$).

Recall first the definitions of internality, analysability and orthogonality.
\defn\label{defn} Let $\pi$ be a partial type over $A$. Then $\pi$ is\begin{itemize}
\item ({\em almost}) {\em $\Sigma$-internal} if for every realization $a$
of $\pi$ there is $B\ind_Aa$ and a tuple $\bar b$ of realizations of types in $\Sigma$
based on $B$, such that $a\in\dcl(B\bar b)$ (or $a\in\bdd(B\bar b)$,
respectively).
\item {\em $\Sigma$-analysable} if for any realization $a$ of $\pi$ there are
$(a_i:i<\alpha)\in\dcl(Aa)$ such that $\tp(a_i/A,a_j:j<i)$ is
$\Sigma$-internal for all $i<\alpha$, and
$a\in\bdd(A,a_i:i<\alpha)$.\end{itemize}
Finally, $p\in S(A)$ is {\em orthogonal} to $q\in S(B)$ if for all $C\supseteq AB$, $a\models p$, and $b\models q$ with $a\ind_A C$ and $b\ind_B C$ we have $a\ind_C b$. The type $p$ is {\em orthogonal to $B$} if it is orthogonal to all types over $B$.
\edefn

We shall say that $a$ is (almost) $\Sigma$-internal or $\Sigma$-analysable over $b$ if $\tp(a/b)$ is.

Now we introduce a general version of the canonical base property.

\defn A simple theory has the {\em canonical base property} with respect to $\Sigma$, if (possibly after naming some constants) whenever $\tp(\cb(a/b))$ is $\Sigma$-analysable then $\tp(\cb(a/b)/a)$ is almost $\Sigma$-internal.\edefn

\bem When $\Sigma$ is the family of types of $SU$-rank $1$ and $a,b$ have finite $SU$-rank, this corresponds to the usual canonical base property. Indeed, since any type of finite $SU$-rank is analysable in the family of types of $SU$-rank $1$, the hypothesis of $\Sigma$-analysability will always be satisfied.\ebem

In \cite{bm-pw} a similar property, named $1$-tight, is defined, but without the condition on $\tp(\cb(a/b))$ to be $\Sigma$-analysable. Instead, $1$-tightness is defined for a family of partial types. Our canonical base property with respect to $\Sigma$ is equivalent to the family of all $\Sigma$-analysable types being $1$-tight with respect to $\Sigma$.

Recall that if $\Sigma$ consists of partial types over $\emptyset$ in a stable theory, if $p(x)\in S(\emptyset)$ is $\Sigma$-internal, there is an $\emptyset$-type-definable group $G$ and a faithfully $\emptyset$-definable action of $G$ on the set of realizations of $p$ which is isomorphic (as a group action) to the group of permutations induced by the automorphisms of $\mathfrak M$ fixing $\Sigma$ pointwise \cite[Th\'eor\`eme 2.26]{udi90, Po87}. We call such a group the {\em Galois} group of $p$ with respect to $\Sigma$.

In a stable theory an $A$-type-definable group is said to be {\em rigid} if all its type-definable connected subgroups are type-definable over $\acl(A)$. Hrushovski, Pillay and the first author have obtained a strong version of the CBP for certain families of types when all Galois groups are rigid.

\tats\label{factSL}\cite[Theorem 2.5]{hpp} Let $\Sigma$ be a family of partial types over $\emptyset$ in a stable theory, and assume all Galois groups (with respect to $\Sigma$) are rigid. If $\tp(\cb(a/b))$ is $\Sigma$-analysable then $\cb(a/b)\subseteq\acl(a,\Sigma)$.\etats

\lmm\label{indep} Let $X$ be an $\emptyset$-invariant set. Then for any $a$ we have
$$a\ind_{\bdd(a)\cap\dcl(X)}X.$$\elmm
\bew Consider a small subset $X_0$ of $X$ such that $a\ind_{X_0} X$ and put $a_0=\cb(X_0/a)$.
Then $a_0$ is definable over a Morley sequence $I$ in $\lstp(X_0/a)$, and $I\subset X$ by invariance of $X$. Therefore $a_0\in\dcl(X)$, whence $a_0\in\bdd(a)\cap\dcl(X)$. As $a\ind_{a_0}X$ by transitivity, the conclusion follows.\qed

\satzli Let $\Sigma$ be a family of partial types over $\emptyset$ in a simple theory. The following are equivalent:
\begin{enumerate} \item Whenever $\tp(b/a)$ is $\Sigma$-analysable, $\cb(a/b)\in\bdd(a,\Sigma)$.
\item Whenever $\tp(b/a)$ is $\Sigma$-analysable, if $c\in\bdd(a,\cb(a/b))$ is $\Sigma$-internal over $a$, then either $c\in\bdd(a)$ or $c\nind_a \Sigma$.
\end{enumerate}
If the theory is stable, both conditions are equivalent to:
\begin{enumerate}
\item[(3)] Whenever $\tp(b/a)$ is $\Sigma$-analysable, if $c\in\acl(a,\cb(a/b))$ is $\Sigma$-internal over $a$, then the connected component of the Galois group of $\stp(c/a)$ with respect to $\Sigma$ acts trivially.
\end{enumerate}
\esatzli
\bew $(1)\Rightarrow(2)$ is immediate. For the other direction, assume $(2)$ and consider some $\Sigma$-analysable $b$ over $a$. Put $$a_0=\bdd(a,\cb(a/b))\cap\bdd(a,\Sigma).$$
Then $\cb(a/b)\ind_{a_0}\Sigma$ by Lemma \ref{indep}, working over $a_0$. Assume $\cb(a/b)\notin\bdd(a_0)$. As $b$ is $\Sigma$-analysable over $a_0$, so is $\cb(a/b)$. Hence there is some $\Sigma$-internal $c\in\bdd(a_0,\cb(a/b))\setminus\bdd(a_0)$ over $a_0$; note that $c\in\bdd(a_0,\cb(a_0/b))$ since $a\subseteq a_0$. As $c\notin\bdd(a_0)$, we have $c\nind_{a_0}\Sigma$ by hypothesis, whence $\cb(a/b)\nind_{a_0}\Sigma$, a contradiction.

Finally, if $T$ is stable and $\tp(c/a)$ is $\Sigma$-internal, then the Galois group of $\stp(c/a)$ with respect to $\Sigma$ acts transitively. Thus its connected component acts trivially if and only if $c\in\acl(a,\Sigma)$. Therefore $(1)\Rightarrow(3)\Rightarrow(2)$.
\qed

Recall that a stable (or simple) theory is {\em non-multidimensional} if every type is non-orthogonal to $\emptyset$. Hence there is a set $A$ of parameters such that every type of finite $SU$-rank is non-orthogonal to the family $\Sigma$ of types over $A$ of $SU$-rank $1$. Thus Fact \ref{factSL} yields immediately the CBP for non-multidimensional stable theories with rigid Galois groups, in the strong form that $\cb(a/b)\subseteq\acl(a,A,\Sigma)$ whenever $SU(\cb(a/b)/A)$ is finite.

Using the following result, we can drop the non-multidimensionality hypothesis.

\tats\cite[Theorem 1.16]{zoe} If $a$ and $b$ are boundedly closed in a simple theory and $SU(\cb(a/b)/a\cap b)$ is finite, then there are $b_1,\ldots,b_m$ independent over $a\cap b$ and types $p_1,\ldots,p_m$ of $SU$-rank $1$, such that $\tp(b_i/a\cap b)$ is analysable in the family of conjugates of $p_i$ over $a\cap b$ for all $i\le m$ and $\bdd(\cb(a/b))=\bdd(b_1,\ldots,b_m)$.\etats

In fact, Chatzidakis formulates her result with algebraic instead of bounded closure, as she assumes elimination of hyperimaginaries (but the generalisation is immediate). Since $p_i$ is non-orthogonal to $\tp(b_i/a\cap b)$, it follows that
$\tp(\cb(a/b)/a\cap b)$ is analysable in the family of non one-based types of $SU$-rank $1$ non-orthogonal to $a\cap b$.

Now if $T$ is stable and all Galois groups with respect to families of types of $SU$-rank $1$ are rigid, then given a canonical base $\cb(a/b)$ of finite $SU$-rank, we can choose a set $A$ of parameters independent of $ab$ over $\acl(a)\cap\acl(b)$ such that every type of $SU$-rank $1$ non-orthogonal to $\acl(a)\cap\acl(b)$ is non-orthogonal to the family $\Sigma$ of types of $SU$-rank $1$ over $A$. So $\cb(a/b)\subseteq\acl(a,A,\Sigma)$ by Fact \ref{factSL}. In particular $\tp(\cb(a/b)/a)$ is almost $\Sigma$-internal, i.e.,\ $T$ has the CBP.

\frag Is the result true for simple theories? \efrag

The problem with a generalization to simple theories is that the Galois group given in \cite{byw} is only {\em almost} hyperdefinable.

\section{Stabilizers and rigidity}
From now on, $G$ will be an $\emptyset$-hyperdefinable group in a simple theory $T$. We shall consider an $A$-hyperdefinable subgroup $H$ of $G$ and an element $g\in G$.

The {\em canonical parameter} $g_H$ of the coset $gH$ over $A$ is the class of $g$ modulo the $A$-hyperdefinable equivalence relation given by $x^{-1}y\in H$. Similarly, we define the canonical parameter for the right coset. Even though the canonical parameter $g_H$ of $gH$ is an $A$-hyperimaginary, there is a hyperimaginary which is interdefinable with $g_H$ over $A$, see \cite[Remark 3.1.5]{wa00}. Working over $A$ we may thus identify the canonical parameter $g_H$ with an ordinary hyperimaginary.

However, in general it need not be true that a hyperdefinable subgroup has a canonical parameter, since equality of hyperdefinable subgroups need not be a type-definable equivalence relation on their parameters.
For canonical parameters to exist we have to assume {\em local connectivity}.

Recall that two subgroups $H_1$ and $H_2$ of $G$ are {\em commensurable} if their intersection has bounded index both in $H_1$ and in $H_2$.
A subgroup $H$ of $G$ is {\em locally connected} if for any model-theoretic or group-theoretic conjugate $H^*$ of $H$, either $H=H^*$ or $H\cap H^*$ has unbounded index in $H$. For every hyperdefinable subgroup $H$ of $G$ there exists a unique minimal hyperdefinable locally connected subgroup commensurable with $H$, its {\em locally connected component} $H^{lc}$, see \cite[Corollary 4.5.16]{wa00}. An inspection of the proof yields that such a subgroup is hyperdefinable over the parameters needed for $H$. Moreover, a locally connected hyperdefinable subgroup, or a coset thereof, has a canonical parameter \cite[Lemma 4.5.19]{wa00}.

\defn Put
$$S(g_H/A)=\{h\in G:\exists x\,[x_H\ind_Ah\land x_H\equiv_A^{Lstp} (hx)_H\equiv_A^{Lstp}g_H]\},$$
and $\stab(g_H/A)=S(g_H/A)\cdot S(g_H/A)$, the {\em (left) stabilizer of $g_H$ in $G$}.\edefn
\satzli $S(g_H/A)$ is hyperdefinable over $\bdd(A)$, and $\stab(g_H/A)$ is a hypderdefinable subgroup of $G$ whose generic types are contained in $S(g_H/A)$.\esatzli
\bew This is an immediate adaptation of \cite[Lemma 4.5.2 and Corollary 4.5.3]{wa00}.\qed

\bem We also have
$$S(g_H/A)=\{h\in G: h\pi\cup\pi\mbox{ does not fork over $A$}\},$$ where $\pi(x)$ is given by
$$\pi(x)=\lstp(g_H/A)=\exists y\,[y^{-1}x\in H\land y\models\lstp(g/A)].$$\ebem

\lmm\label{stabilizers} Let $(H_i:i<\alpha)$ be a continuous descending sequence of $A$-hyperdefinable subgroups of $G$, i.e.\ $H_i=\bigcap_{j<i}H_{j+1}$ for $i<\alpha$. If $g\in G$, then\begin{enumerate}
\item\label{stab1} for $i\le j$ we have $\stab(g_{H_i}/A)\ge\stab(g_{H_j}/A)$, and
\item\label{stab2} for $\lambda$ limit, $\stab(g_{H_\lambda}/A)=\bigcap_{i<\lambda}\stab(g_{H_i}/A)$.\end{enumerate}\elmm
\bew Suppose $h\ind_A g_{H_j}$ with $(hg)_{H_j}\equiv_A^{Lstp}g_{H_j}$. Since $H_i\ge H_j$ we have $h\ind_A g_{H_i}$ and $(hg)_{H_i}\equiv_A^{Lstp}g_{H_i}$. Hence $\stab(g_{H_j}/A)\le\stab(g_{H_i}/A)$.

From the first statement we obtain $\stab(gH_\lambda/A)\le\bigcap_{i<\lambda}\stab(g_{H_i}/A)$. Now suppose $h\in\bigcap_{i<\lambda}\stab(g_{H_i}/A)$ is generic. Then every formula $\phi(x)\in\tp(h/A)$ is generic in $\stab(g_{H_i}/A)$ for $i<\lambda$ sufficiently big. By the first paragraph, there is $h_i\models\phi$ such that $h_i\ind_A g_{H_j}$ and $(h_ig)_{H_j}\equiv_A^{Lstp}g_{H_j}$ for all $j\le i$. By compactness, there is $h'\models\tp(h/A)$ such that $h'\ind_A g_{H_i}$ and $(h'g)_{H_i}\equiv_A^{Lstp}g_{H_i}$ for all $i<\lambda$. But this implies $h'\ind_Ag_{H_\lambda}$ and $(hg)_{H_\lambda}\equiv_A^{Lstp}g_{H_\lambda}$, whence $h'\in\stab(g_{H_\lambda}/A)$.\qed

Recall that for any set $A$, the {\em $A$-connected component} $G_A^0$ of $G$ is the smallest $A$-hyperdefinable subgroup of bounded index; note that it is normal. Whereas in a stable theory the $A$-connected component of a group does not depend on $A$ and is locally connected, this need not be true in a simple theory.

\lmm\label{kernel} If $H$ is normalized by $G^0_A$, then $H^{g^{-1}}\le\stab(g_H/A)$.\elmm
\bew Let $h\in H^{g^{-1}}$ and consider $x\models\lstp(g/A)$ with $x\ind_A h$. Then $xG^0_A=gG^0_A$, as this coset is $\bdd(A)$-hyperdefinable. Hence $h\in H^{g^{-1}}=H^{x^{-1}}$, and $hx\in xH$. Thus $(hx)_H=x_H\equiv_A g_H$; as clearly $x_H\ind_A h$ we get $h\in\stab(g_H/A)$.\qed

\defn\label{def:rigid} An $A$-hyperdefinable group is {\em rigid} if every hyperdefinable subgroup is commensurable with one hyperdefinable over $\bdd(A)$. \edefn

\bem An $A$-hyperdefinable group is rigid if and only if every locally connected subgroup is hyperdefinable over $\bdd(A)$. \ebem

In the pure theory of an algebraically closed field, semi-abelian varieties are rigid. On the other hand, one-based groups are rigid \cite{hp}.

\lmm\label{normal} If $G$ is rigid, then every hyperdefinable locally connected subgroup is normalized by $G^0_\emptyset$.\elmm
\bew Let $H$ be a hyperdefinable locally connected subgroup of $G$. Then every $G$-conjugate of $H$ is also locally connected, and hence hyperdefinable over $\bdd(\emptyset)$. So there are only boundedly many $G$-conjugates, and $G^0_\emptyset$ must normalize them all.\qed

The next two results give some basic properties of rigid groups; they will not be used in the sequel. Proposition \ref{nilpotent} has been noted for groups of finite Morley rank in \cite[Remark 1.12]{hpp}.

\satzli\label{nilpotent} A type-definable superstable rigid group $G$ is nilpotent-by-finite.\esatzli
\bew Clearly, we may assume that $G$ is connected. If $U(G)=\omega^\alpha\cdot n+\beta$ with $\beta<\omega^\alpha$, then by \cite{BL} there is a type-definable connected abelian subgroup $A$ with $U(A)\ge\omega^\alpha$. So $A$ is normal by Lemma \ref{normal}, and $U(G/A)<U(G)$. By inductive hypothesis $G/A$ is nilpotent.

If $g\in G$ is generic, let $C$ be the centralizer-connected component of $C_G(g)$. Then $C$ is locally connected and of finite index $n$ in $C_G(g)$; since it is $\bdd(\emptyset)$-definable, it does not depend on $g$. Therefore $C$ is centralized by all generic elements, and thus by the whole of $G$. Put $N=A\cdot Z(G)$. Then $N$ is normal and abelian, and $G/N$ is nilpotent of exponent at most $n$. So $G$ is nilpotent by a theorem of Baudisch and Wilson \cite{BW}.\qed

Rigidity implies a monotonicity property for stabilizers.

\lmm\label{finite} Suppose $G$ is rigid, $g\in G$ and $A\subseteq B$. Then $\stab(g/B)_B^0$ is contained in $\stab(g/A)$.\elmm
\bew By rigidity $\stab(g/B)_B^0$ is commensurate with a group $H$  hyperdefinable over $\bdd(\emptyset)$; by $B$-connectivity it has bounded index in $H$. Consider a generic $h\in\stab(g/B)_B^0$. Then $h$ is generic in $H$ and so $h\ind B$. Also, since $h$ is generic in $\stab(g/B)$, there is $g'\models\lstp(g/B)$ with $g'\ind_B h$ and $hg'\models\lstp(g/B)$. In particular $g',hg'\models\lstp(g/A)$. Since $h\ind B$ we have $h\ind_A B$, whence $h\ind_A g'$ by transitivity and hence $h\in\stab(g/A)$.\qed

Recall that $\ell_1^\Sigma(g/A)$ is the set of all elements in $\bdd(g,A)$ whose type over $A$ is almost $\Sigma$-internal.

\lmm\label{levels} If $a\ind_A b$ and $\tp(c/A,b)$ is almost $\Sigma$-internal, then $a\ind_B b,c$ where $B=\ell_1^\Sigma(a/A)$.\elmm
\bew \cite[Corollary 2.4]{pw} yields that $\tp(\cb(b,c/A,a)/A)$ is almost $\Sigma$-internal. Thus $\cb(b,c/A,a)\in B$.\qed

\kor Let $g\in G$ and suppose that there is some $h\in G$ such that $gh\ind_A g$ and $\tp(h/A)$ is almost $\Sigma$-internal. Then $g$ is almost $\Sigma$-internal over $A$.\ekor
\bew Put $B=\ell_1^\Sigma(g/A)$. By Lemma \ref{levels} we have $g\ind_B gh,h$, so $g\in\bdd(B)=B$.\qed

In order to obtain equality in Lemma \ref{stabilizers}(\ref{stab1}), we have to assume rigidity, $\Sigma$-analysablility and work with $\ell_1^\Sigma(g/A)$-connected components.

\satzli\label{lem:stab+rigid} Let $G$ be rigid. If $N\trianglelefteq H$ are $A$-hyperdefinable subgroups normalized by $G^0_A$ such that $H/N$ is $\Sigma$-analysable, $g\in G$ and $B=\ell_1^\Sigma(g/A)$, then $\stab(g_H/B)_B^0=\stab(g_N/B)_B^0(H^{g^{-1}})_B^0$.\esatzli
\bew Since $H/N$ is $\Sigma$-analysable, by \cite[Corollary 4.6.5]{wa00} there is a continuous descending sequence $(H_i:i\le\alpha)$ of hyperdefinable normal subgroups of $H$ with $H_0=H$ and $H_\alpha=N$, such that $H_i/H_{i+1}$ is almost $\Sigma$-internal for all $i<\alpha$. As $\Sigma$ is $\emptyset$-invariant and both $H$ and $N$ are normalized by $G_A^0$, we can successively choose the $H_i$ to be hyperdefinable over $A$ and normalized by $G^0_A$.
\beh $\stab(g_{H_i}/B)_B^0=\stab(g_{H_j}/B)_B^0(H_i^{g^{-1}})_B^0$ for all $j\ge i$.\ebeh
\bewbeh Consider a generic element $h\in\stab(g_{H_i}/B)^0_B$. We may assume that $h\ind_Bg$ and $(hg)_{H_i}\equiv_B g_{H_i}$. As in the proof of Lemma \ref{finite}, rigidity implies $h\ind B$, whence $h\ind_A g$ by transitivity.

Choose some $h'\in H_i^{g^{-1}}$ with $hh'g\equiv_B g$. As $G_A^0$ normalizes $H_i$, the quotient $H_i^{g^{-1}}/H_{i+1}^{g^{-1}}$ is hyperdefinable over $\bdd(A)$ and almost $\Sigma$-internal. Thus $\tp(h'_{H_{i+1}^{g^{-1}}}/A)$ is almost $\Sigma$-internal. Hence $h,h'_{H_{i+1}^{g^{-1}}}\ind_Bg$ by Lemma \ref{levels}. Therefore there is some $h''\in hh'H_{i+1}^{g^{-1}}$ with $h''\ind_Bg$. But then
$$h''gH_{i+1}=h''H_{i+1}^{g^{-1}}g=hh'H_{i+1}^{g^{-1}}g=hh'gH_{i+1}.$$
Thus $(h''g)_{H_{i+1}}=(hh'g)_{H_{i+1}}\equiv_Bg_{H_{i+1}}$ and $h''\in\stab(g_{H_{i+1}}/B)$; as $h'\in H_i^{g^{-1}}$ and $h''H_{i+1}^{g^{-1}}=hh'H_{i+1}^{g^{-1}}$ we obtain
$h\in h''H_i^{g^{-1}}$, whence $$h\in\stab(g_{H_{i+1}}/B)\,H_i^{g^{-1}}.$$
Taking $B$-connected components yields
$$\stab(g_{H_i}/B)_B^0\le\stab(g_{H_{i+1}}/B)_B^0(H_i^{g^{-1}})_B^0\ ;$$
the other inclusion follows from Lemmas \ref{stabilizers} and \ref{kernel} by taking $B$-connected components. The claim now follows by transfinite induction from Lemma \ref{stabilizers}(\ref{stab2}).\qed

Taking $i=0$ and $j=\alpha$, we obtain the proposition.\qed

\section{The Socle Lemma}
Hrushosvki's original {\em Socle Lemma}, formulated first for finite Morley rank \cite[Proposition 4.3]{udi} and then for finite $SU$-rank \cite[Proposition 3.6.2]{udi2}, states the following:
\begin{SL}Let $G$ be an abelian group of finite $SU$-rank, $H$ a definable subgroup and $X$ a (type-)definable subset of $G$. Assume:\begin{enumerate}
\item Every $\acl(G/H)$-definable subgroup of $H$ is commensurable with an $\acl(\emptyset)$-definable one ($G/H$-rigidity).
\item $G$ has no definable subgroup $H'\ge H$ with $H'/H$ infinite and $H'\subseteq\acl(H,Y,C)$ for some $SU$-rank one set $Y$ and finite $C$.
\item For any $b\in X$ and parameter set $B$, $\stab(b/B)\cap H$ is finite.\end{enumerate}
Then up to smaller rank, $X$ is contained in finitely many cosets of $H$.\end{SL}
In particular, we deduce:
\begin{SC}Let $G$ be a rigid abelian group of finite $SU$-rank, and $p$ a Lascar strong $G$-type of finite stabilizer. Then $p$ is almost internal in the set $\Sigma$ of all types of $SU$-rank one.\end{SC}
\bew By the Indecomposability Theorem \cite[Theorem 5.4.5]{wa00} there is an almost $\Sigma$-internal definable subgroup $H$ such that any type in $\Sigma$ is contained in finitely many cosets of $H$. This yields (2); hypothesis (1) follows from rigidity, and hypothesis (3) holds by assumption for $X=p$. Hence some non-forking extension of $p$ is contained in finitely many cosets of $H$, and thus is almost $\Sigma$-internal.\qed

As pointed out in the introduction, the above Corollary holds with rigidity replaced by the canonical base property \cite{pz,hpp}. More generally, Chatzidakis proves the following two facts (she works in a simple theory with elimination of hyperimaginaries and the family $\Sigma$ of all non locally modular (Lascar) strong types of $SU$-rank $1$, but the generalization is straightforward).

\tats\label{fact:socle}\cite[Proposition 2.6]{zoe} Suppose $T$ has the CBP with respect to $\Sigma$ and $G$ is $\Sigma$-analysable. If $g\in G$ and $\bar g$ is the canonical parameter of $\stab(g/A)\, g$, then $\tp(\bar g/A)$ is almost $\Sigma$-internal.\etats

\tats\cite[Corollary 2.7]{zoe} Suppose $T$ has the CBP with respect to $\Sigma$ and $G$ is $\Sigma$-analysable. If $g\in G$ then $\stab(g/\ell_1^\Sigma(g/A))\le \stab(g/A)$.\etats

When $\stab(g/A)$ is bounded, Fact \ref{fact:socle} yields the CBP version of the Corollary above, see for example \cite[Fact 1.3]{hpp}. We shall now prove a rigidity version in the spirit of the original Socle Lemma.

\satz\label{thm:generic} Suppose $G$ is $\Sigma$-analysable and rigid. If $g\in G$ and $B=\ell_1^\Sigma(g/A)$, then $\stab(g/B)\,g$ is $B$-hyperdefinable.\esatz
\bew Since $G$ is $\Sigma$-analysable, there is a continuous descending sequence $(G_i:i\le\alpha)$ of $\emptyset$-hyperdefinable normal subgroups of $G$ with $G_0=G$ and $G_\alpha=\{1\}$ such that $G_i/G_{i+1}$ is almost $\Sigma$-internal for all $i<\alpha$. Proposition \ref{lem:stab+rigid} with $H=G_i$ and $N=\{1\}$ yields for $i<\alpha$ that
$$\stab(g_{G_i}/B)_B^0=\stab(g/B)_B^0(G_i)_B^0.$$
We shall show inductively that $\stab(g_{G_i}/B)\,g$ is hyperdefinable over $B$. For $i=\alpha$ this yields the result.

The assertion is clear for $i=0$; for limit ordinals it follows from continuity of the sequence $(G_i:i\le\alpha)$ and Lemma \ref{stabilizers}(\ref{stab2}). So assume inductively that the coset $\stab(g_{G_i}/B)\,g$ is hyperdefinable over $B$.
Then also $\stab(g_{G_i}/B)^0_B\,g$ is hyperdefinable over $B$, as $B$ is boundely closed. If $g'\models\tp(g/B)$, then
$$g'g^{-1}\in\stab(g_{G_i}/B)_B^0\le\stab(g_{G_i}/B)^{lc}.$$
By rigidity both $\stab(g_{G_i}/B)^{lc}$ and $\stab(g_{G_{i+1}}/B)^{lc}$ are hyperdefinable and connected over $\bdd(\emptyset)$ and normalized by $G_\emptyset^0$. In particular,
$$\stab(g_{G_{i+1}}/B)^{lc}$$
is normalized by $\stab(g_{G_i}/B)^{lc}$. Now
$$Q=\stab(g_{G_i}/B)^{lc}/(\stab(g_{G_{i+1}}/B)^{lc}\cap\stab(g_{G_i}/B)^{lc})$$
is $\bdd(\emptyset)$-hyperdefinable. Moreover, $Q$ is isogenous to
$$\begin{aligned}\stab(g_{G_i}/B)^0_B&/(\stab(g_{G_{i+1}}/B)^0_B\cap\stab(g_{G_i}/B)^0_B)\\
&=\stab(g/B)_B^0(G_i)_B^0/\stab(g/B)_B^0(G_{i+1})_B^0\end{aligned}$$
(note that in the hyperdefinable context, an isogeny allows for {\em bounded} rather than {\em finite} indices and kernels).
It is thus a homomorphic image of $(G_i)^0_B/(G_{i+1})^0_B$; as $G_i/G_{i+1}$ is almost $\Sigma$-internal, so is $Q$, and hence also $\tp((g'g^{-1})_{(\stab(g_{G_{i+1}}/B)^{lc}\cap\stab(g_{G_{i}}/B)^{lc})})$.
In particular,
$$\tp((g'g^{-1})_{\stab(g_{G_{i+1}}/B)^{lc}}/A)$$ is almost $\Sigma$-internal, whence by Lemma \ref{levels}
$$(g'g^{-1})_{\stab(g_{G_{i+1}}/B)^{lc}}\ind_Bg.$$
Put $H=\stab(g_{G_{i+1}}/B)^{lc}$. Then $H$ and $H^g$ are both $B$-hyperdefinable, and
$$(g'g^{-1})H\,gH^g=g'H^g.$$
It follows that
$$(g'g^{-1})H\le\stab(g_{H^g}/B).$$
In fact, since $g$ and $g'$ have the same Lascar strong type over $B$, we even have
$$(g'g^{-1})H^0_B\le\stab(g_{H^g}/B)^0_B.$$
By Proposition \ref{lem:stab+rigid} with $N=\{1\}$ we obtain
$$\stab(g_{H^g}/B)^0_B=\stab(g/B)^0_B\,H^0_B=H^0_B,$$
where the second equality follows from Lemma \ref{stabilizers}(\ref{stab1}). Hence
$$(g'g^{-1})H^0_B=H^0_B\quad\text{and}\quad g^{-1}H^0_B=g'^{-1}H^0_B.$$
Thus $g^{-1}H^0_B$ is hyperdefinable over $B$ , as is $g^{-1}H=g^{-1}H^0_B\,H$, and also $Hg=(g^{-1}H)^{-1}$, and finally
$\stab(g_{G_{i+1}}/B)\,g=\stab(g_{G_{i+1}}/B)\, Hg$.\qed

\bem In fact, hyperdefinability of the coset implies genericity:
If $\stab(g/A)\,g$ is $A$-hyperdefinable, then $tp(g/A)$ is generic in it.\ebem
\bew Suppose $\stab(g/A)\,g$ is $A$-hyperdefinable. Then for any stratified local rank $D$ we have
$$D(\stab(g/A)\,g)\ge D(g/A)\ge D(\stab(g/A))=D(\stab(g/A)\,g).$$
Thus equality holds, and $g$ is generic in the coset $\stab(g/A)\,g$.\qed

\kor\label{socle} Let $G$ be an $\emptyset$-hyperdefinable group and $g\in G$. If $G$ is  $\Sigma$-analysable and rigid, then whenever $\stab(g/A)$ is bounded, $\tp(g/A)$ is almost $\Sigma$-internal.\ekor
\bew By Lemma \ref{stabilizers} $\stab(g/B)$ is bounded and hence $g$ is bounded over $\ell_1^\Sigma(g/A)$ by Theorem \ref{thm:generic}, i.e.\ $\tp(g/A)$ is almost $\Sigma$-internal.\qed

\bem It follows from the proof of Theorem \ref{stabilizers} that we can weaken the rigidity hypothesis, either to $\stab(g_{G_i}/B)/G_i$ being bounded for all $i\ge2$, or to $\stab(g/B)$ being bounded and $\stab(g_{G_i}/B)$ being commensurable to a group hyperdefinable over $A$ whenever $2\le i<\alpha$.\ebem

In particular, Theorem \ref{thm:generic} holds without any rigidity hypothesis when $G$ is $\Sigma$-analysable in two steps:

\bem Let $G$ be an $\emptyset$-hyperdefinable group, $g\in G$ and $B=\ell_1^\Sigma(g/A)$. If $G$ is $\Sigma$-analysable in two steps, then $\tp(g/B)$ is generic in a $B$-hyperdefinable coset of its stabilizer.\ebem
\bew We sketch a short proof for convenience. If $G$ is $\Sigma$-analysable in two steps, then there is some hyperdefinable normal subgroup $N$ of $G$ such that $G/N$ and $N$ are almost $\Sigma$-internal. Thus $gN\in B$, and for any realization $g'\models\tp(g/B)$ we have $g'g^{-1}\in N$. Hence $\tp(g'g^{-1}/A)$ is almost $\Sigma$-internal, and $g\ind_B g'g^{-1}$ by Lemma \ref{levels}. Thus $g'g^{-1}\in\stab(g/B)$. Therefore $\stab(g/B)\,g$ is hyperdefinable over $B$ and the type $\tp(g/B)$ must be generic in it. \qed

\bem In the finite rank context, we can extend the Corollary of the Socle Lemma to non-abelian rigid groups: By the Indecomposability Theorem there is a normal almost $\Sigma$-internal $\emptyset$-hyperdefinable subgroup $N$ such that any almost $\Sigma$-internal partial type $\pi$ is contained in only finitely many cosets of $N$. Then a Lascar strong type is almost $\Sigma$-internal if and only if it is contained in a single coset of $N$; by Corollary \ref{socle} this holds in particular for any Lascar strong type with bounded stabilizer.\ebem

\section{Final remarks}
As we have seen, there are two ways to deduce almost internality of types with finite stabilizer from rigidity: Either by assuming that the Galois groups are rigid, deducing the canonical base property \cite[Theorem 2.1]{hpp} and using Fact \ref{fact:socle}, or by assuming directly that the ambient group is rigid (Corollay \ref{socle}). A priori neither hypothesis implies the other one. It would thus be interesting to compare the two approaches, and to identify connections beween a group $G$ and the Galois groups of types in $G$.

Another question concerns nilpotent groups: In \cite[Remark 6.7]{pw} it is shown that if $G$ is $\Sigma$-analysable type-definable or supersimple, then there is a nilpotent normal subgroup $N$ such that $G/N$ is almost $\Sigma$-internal. Is this related to Proposition \ref{nilpotent} stating that a superstable rigid group is virtually nilpotent?

\end{document}